\documentclass[12pt]{article}
\usepackage{amssymb,amsbsy,amsmath,amsfonts,amssymb}
\usepackage{latexsym,euscript,exscale}

\usepackage{times}

\newcommand{\pf}{

\smallskip

\noindent {\it Proof : }}

\newcommand {\N}{\mathbb N}

\newcommand {\R}{\mathbb R}

\newcommand{\pff}{$\hfill \square$
\smallskip}

\newcommand{\norm}[1]{\ensuremath{\left\|#1\right\|}}

\newtheorem{prop}{Proposition}
\newtheorem{lemm}[prop]{Lemma}
\newtheorem{theo}[prop]{Theorem}

\newtheorem{defi}[prop]{Definition}

\title{A Banach space dichotomy theorem for quotients of subspaces}

\author{Valentin Ferenczi}

\date{}

\begin{document}

\maketitle

\begin{abstract}
A Banach space $X$ with a Schauder basis is defined to have the {\em restricted
quotient here\-ditarily indecomposable
 property} if 
$X/Y$ is hereditarily indecomposable for any infinite codimensional subspace $Y$ with a successive
finite-dimen\-sio\-nal decomposition on
the basis of $X$. The following dichotomy theorem is proved: any infinite dimensional Banach space contains a quotient of subspace which either
has an unconditional basis, or has the restricted quotient hereditarily indecomposable
property. 
\end{abstract}

\section{Introduction}

In 2002, W.T. Gowers published his famous Ramsey theorem for block-subspaces in a Banach space \cite{G1}.
If $X$ is a Banach space with a Schauder basis, {\em block-vectors} in $X$ denote
 non zero vectors with finite support on the basis, and
 {\em block-sequences} are infinite sequences of block-vectors with successive supports;  {\em block-subspaces} are subspaces
generated by block-sequences.

If $Y$ is a block-subspace of $X$, {\em Gowers' game in $Y$} is the infinite game where
Player 1 plays  block-subspaces $Y_n$ of $Y$, and Player 2  plays normalized block vectors
$y_n$ in $Y_n$. 

If $\Delta=(\delta_n)_{n \in \N}$ is a sequence of reals,    $\Delta>0$ means that
$\delta_n>0$ for all $n \in \N$.
For $A$  a set of normalized block-sequences, and any $\Delta=(\delta_n)_{n \inÊ\N}>0$,
 let $A_{\Delta}$ be the set of normalized block-sequences $(y_n)_{n \in \N}$ such that
there exists $(x_n)_{n \in \N}$ in $A$ with $\norm{x_n-y_n} \leq \delta_n$ for all $n \in \N$.

\begin{theo}(Gowers' Ramsey Theorem) Let $A$ be a set of normalized block-sequences which is analytic as a subset of
$X^{\omega}$ with the product of the norm topology on $X$. Assume that every block-subspace of $X$
contains a block-sequence in $A$. Let $\Delta>0$.
Then there exists a block-subspace $Y$ of $X$ such that 
Player 2 has a winning strategy in Gowers' game in $Y$ for producing a sequence $(y_n)_{n \in \N}$ in $A_{\Delta}$.
\end{theo}

The most important consequence of the Ramsey
 Theorem of Gowers is the
so-called dichotomy theorem for Banach spaces.
A Banach space $X$ is said to be {\em decomposable} if it is 
 a direct (topological) sum of two infinite-dimensional closed subspaces.
An infinite dimensional space is  {\em hereditarily indecomposable (or HI)} 
when it has no decomposable subspace. A Schauder basis $(e_n)_{n \in \N}$ of $X$ is 
unconditional if there exists $C \geq 1$ such that for all
 $\sum_{i \in \N} \lambda_i e_i$ in $X$, all
$(\epsilon_i)_{i \in \N} \in \{-1,1\}^\N$, $\norm{\sum_{i \in \N} \epsilon_i \lambda_i e_i} \leq
C\norm{\sum_{i \in \N} \lambda_i e_i}$. 

\begin{theo}(Gowers' Dichotomy Theorem)
Every infinite dimensional Banach space contains a subspace $Y$ which satisfies one of the two 
following properties, which are both possible, and mutually exclusive:

i) $Y$ has an unconditional basis,

ii) $Y$ is hereditarily indecomposable. 
\end{theo}

\

These properties are even exclusive in the sense that if a space satisfies i) (resp. ii)), then no further
subspace satisfies ii) (resp. i)). Indeed if a Banach space $X$ is hereditarily indecomposable,
then so is any subspace
of $X$; and if $X$ has an unconditional basis, then every block-subspace of $X$ has an 
unconditional basis, and so any subspace of $X$ has a further subspace with
an unconditional basis.

\subsection{HI spaces and their quotient spaces}

From now on, spaces and subspaces are supposed infinite dimensional and closed unless specified otherwise.
For two subspaces $Y$ and $Z$ of a space $X$, a convenient notion of angle was used by B. Maurey to give a 
simple proof of Gowers' dichotomy theorem \cite{M2}: let
$$a(Y,Z)=\inf_{y \in Y, z \in Z, y \neq z}\frac{\norm{y-z}}{\norm{y+z}}.$$
It is in particular clear that $a(Y,Z) \neq 0$ if and only if $Y+Z$ forms a topological
 direct sum
 in $X$, and therefore a space $X$ is hereditarily indecomposable if and only if
$a(Y,Z)=0$ for any subspaces $Y,Z$ of $X$. On the other hand, a basic sequence $(e_i)_{i \in \N}$
 is $C$-unconditional
if $a([e_i, i \in I],[e_i, i \in J]) \geq 1/C$ for every partition $\{I,J\}$ of $\N$,  
where $[e_i, i \in I]$ denotes the closed linear space generated by $(e_i)_{i \in I}$.

 We also note that it was proved in \cite{GM1} that hereditarily indecomposable spaces 
are never isomorphic to  proper subspaces.

\

While classical spaces, such as $c_0$ and $\ell_p, 1 \leq p <+\infty$, or $L_p,1<p<+\infty$,
 have  unconditional bases,
the first known 
example of a HI space was given by Gowers and Maurey in 1993, \cite{GM1}. Gowers-Maurey's space
$X_{GM}$ is actually {\em quotient hereditarily indecomposable (or QHI)}, that is,
no quotient of a subspace of $X_{GM}$ is decomposable, or equivalently, every infinite dimensional
quotient
space of $X_{GM}$ is HI \cite{F3}; as $X_{GM}$ is reflexive,
 it follows  that $X^*_{GM}$ is also quotient hereditarily indecomposable, and in particular also 
hereditarily indecomposable.
 In \cite{F3}, an example $X$ was also provided 
which is HI and not QHI. This example is defined as the "push-out"
 $(X_1 \oplus X_2)/\{(y,-y): y \in Y\}$
of two specific Gowers-Maurey's type spaces $X_1$ and $X_2$ with respect to  a "common" subspace
 $Y$. It is therefore still very close to being QHI, in the sense that it is saturated
with QHI subspaces, and the natural quotient space of $X$ which is decomposable is a direct
sum of two HI spaces. This led the author to conjecture that any
 quotient of a HI space should contain a HI or even QHI subspace, or that the dual of any reflexive HI space 
should contain a HI subspace.

\

This however turned out to be completely false.
 Examples of HI spaces were built with quotients which are very far from being
 HI.
 Using methods based on the definition of some notion of HI interpolation of Banach spaces,
 S. Argyros and V. Felouzis constructed a HI space
with some quotient space isomorphic to
 $c_0$ (resp. $l_p, 1<p<+\infty$) \cite{AF}. S. Argyros and A. Tolias used deep constructions, based
on what is now known 
as the "extension
method" \cite{AAT}, to prove that any separable
 Banach space which
does not contain a copy of $l_1$ is isomorphic to the quotient space of some separable HI space
\cite{AT}; and to construct a  reflexive Banach space $X$ which is HI but whose dual is
saturated with unconditional basic sequences \cite{AT2}, therefore any quotient space of $X$
 has a further quotient with an unconditional basis. These results shatter all hopes of
general results preserving the HI property when passing to quotient spaces, or to the dual. 
We refer to \cite{A}, \cite{AT},  and \cite{M} for more details about these examples and hereditarily 
indecomposable spaces in general, as well as about other examples, and also to the recent work \cite{AAT} which contains
an comprehensive  introduction to the previous examples.

S. Argyros
 asked whether there existed a reflexive HI Banach space $X$, such that no subspace of $X$ has a HI dual.
 This would show that the H.I. structure is in general  not inherited by duals, not even 
in a very weak sense.
 None of the 
HI examples constructed so far seem to answer that question
 (for more about this, we refer to the remarks and questions section at the end of this paper). 

\

Our main result is somewhat related to the question of Argyros.
Its starting point is the observation that the
 situation becomes more pleasant again when one looks at quotient of subspaces (or QS-spaces) of a
given Banach space.
First note that the features of the QHI property
 with respect to quotient of subspaces are quite similar to the ones
 of the HI property with respect to subspaces. Indeed this property obviously passes to further
 QS-spaces. We also have the following result.

\begin{prop} If $X$ is hereditarily indecomposable, then $X$ 
is isomorphic to no proper quotient of subspace of itself. \end{prop} 
\pf Assume $X$ is HI and $\alpha$ is an isomorphism from $X$ onto $Y/Z$
for some $Z \subset Y \subset X$. We may assume that $\dim Z=+\infty$.
Then by properties of HI spaces \cite{M2}, the quotient map $\pi:Y \rightarrow Y/Z$ is strictly
singular. The map $T=\alpha\pi$ is an onto map whose Fredholm index $i(T)$ (defined
as $\dim(Ker T)-\dim(X/TY)$ when this expression has a meaning) is $+\infty$. By continuity of
the index (\cite{LT} Proposition 2.c.9), we deduce that $i(T-\epsilon i_{YX})=+\infty$ for some
small enough $\epsilon>0$. On the other hand, $T$ is strictly singular, therefore, by \cite{LT}
Proposition 2.c.10, $i(T-\epsilon i_{YX})=i(-\epsilon i_{YX}) \leq 0$. \pff

\

The unconditional property also satisfies some type of heredity for quotient of subspaces.
T. Odell proved that if $X$ has a shrinking finite-dimensional unconditional decomposition,
 then every normalized weakly null sequence in a quotient of $X$ has an unconditional 
subsequence \cite{O}, and therefore every QS-space of $X$ contains
 an unconditional basic sequence. 

\

It is therefore tempting to look for some general dichotomy result for quotient of subspaces involving the QHI 
property on one side and some unconditionality property on the other.

\

\subsection{Angles between quotient of subspaces}

To motivate our following definitions and
results, we take a closer look at Gowers-Maurey's sequence space $X_{GM}$.
To prove that $X_{GM}$ is HI, Gowers and Maurey
build, for arbitrary large $k \in \N$, successive  biorthonormal sequences $(y_i)_{i \leq k}$ and 
$(y_i^*)_{i \leq k}$ of "special" pairs of
vectors and functionals, such that   
$$\norm{\sum_{i \leq k}y_i^*} \simeq  \sqrt{\log(k)},$$ while
$$\norm{\sum_{i \leq k}(-1)^i y_i} \simeq k/\log(k).$$ 
Up to a perturbation,
 the terms $(y_i)_{i \leq k}$ may be taken in arbitrary subspaces of $X_{GM}$.
Therefore, given $Y,Z \subset X_{GM}$, by taking the  even terms close enough to $Y$ and 
 the odd terms close enough to $Z$,
 we may find
 vectors
$y$ almost in $Y$ and $z$ almost in $Z$, and functionals
$y^*$ and $z^*$, with disjoint supports,
 such that $\norm{y-z} \simeq k/\log(k)$ while
 $$\norm{y+z} 
\geq \frac{(y^*+z^*)(y+z)}{\norm{y^*+z^*}} \simeq k/\sqrt{\log(k)}.$$ It follows that $Y+Z$
 is never a direct sum.

\

The proof in \cite{F3} that $X_{GM}$ is QHI is based on the fact 
 that one can actually choose $y^*$ and $z^*$ close enough to
 $W^{\perp}$ for any $W$ which is an infinite codimensional subspace
 of $Y$ and of $Z$. It follows easily that $X_{GM}$ is quotient hereditarily indecomposable.
By the proof it is clear than one can even pick $y^*$ close enough to $V^{\perp}$ and
$z^*$ close enough to $W^{\perp}$ for any 
  infinite codimensional subspaces $V$ of $Y$ and $W$ of $Z$.  The point here is that each term
  of the sequences of "special"  vectors (resp. functionals) must be taken
 in some  set $A_n$ (resp. $A_n^*$) which is asymptotic, i.e. 
intersects any subspace of $X_{GM}$ (resp. $X^*_{GM}$), associated to some $n$ depending on the previous terms,
 but the subspace 
in which to pick it may be chosen arbitrarily.

\

 For $X$ a Banach space, and a subspace $Y_*$ of $X^*$, denote by $\norm{.}_{Y_*}$ the semi-norm
defined on $X$ by $\norm{x}_{Y_*}=\sup_{y^* \in Y_*, \norm{y^*} \leq 1} y^*(x)$, and by
 $Y_*^{\perp}$  the orthogonal of $Y_*$ in $X$. When $Y_*=Y^{\perp}$ for some $Y \subset X$,
 $\norm{.}_{Y_*}$ is the quotient norm on
$X/Y$.

A {\em QS-pair} is some $(Y_*,Y) \subset X^* \times X$ such that $Y_*^{\perp} \subset Y$. It
 may be associated 
to the QS-space $Y/Y_*^{\perp}$. The natural notion of inclusion between QS-pairs
$$(Z_*,Z) \subset (Y_*,Y) \Leftrightarrow (Z_* \subset Z) \wedge (Y_* \subset Y)$$
corresponds to taking quotient of subspaces of the associated QS-spaces.
Indeed if
$(Z_*,Z) \subset (Y_*,Y)$, 
then $Z/Z_*^{\perp} \simeq (Z/Y_*^{\perp})/(Z_*^{\perp}/Y_*^{\perp})$.
An {\em infinite dimensional QS-pair} is a QS-pair whose associated
QS-space is infinite dimensional.
We
define the angle $A((Y_*,Y),(Z_*,Z))$ between two QS-pairs 
by 

$$A((Y_*,Y),(Z_*,Z))=\inf_{y \neq z,y^* \neq z^*,y^*(z)=z^*(y)=0}\frac{\norm{y-z}\norm{y^*-z^*}}
 {|y^*(y)-z^*(z)|},$$ 
where the infimum is taken over $y \in Y, z \in Z, y^* \in Y_*, z^* \in Z_*$.

\

Note that if we let $W_*=Y_*=Z_*$, then we obtain 
$$A((W_*,Y), (W_*,Z)) \geq \inf_{y \neq z,\norm{y^*-z^*}=1}\frac{\norm{y-z}}
 {|(y^*-z^*)(y+z)|} \geq \inf_{y \neq z}\frac{\norm{y-z}_{W_*}}{\norm{y+z}_{W_*}},$$
and therefore when $W_*=W^{\perp}$ for some $W \subset X$,
$A((W_*,Y),(W_*,Z)) \geq  a(Y/W ,Z/W).$
 In particular 
$Y/W$ and $Z/W$ do not form a direct sum in $X/W$ when
$A((W^{\perp},Y),(W^{\perp},Z))=0$. If this is true for all $W,Y,Z$ with $W$ an infinite
codimensional subspace of $Y$ and of $Z$
then we deduce that
$X$ is QHI.

\

By our previous description of special sequences in Gowers-Maurey's space, $X_{GM}$ is an exemple
 of a reflexive space for which $A((Y_*,Y),(Z_*,Z))=0$ for all infinite dimensional
QS-pairs $(Y_*,Y)$ and $(Z_*,Z)$ of $X$. Indeed if $y$, $z$, $y^*$, $-z^*$ are the odd and even 
parts respectively of adequate length $k$ special sequences, we have
$$\norm{y-z}\norm{y^*-z^*} \simeq k/\sqrt{\log(k)},$$
while $$|y^*(y)-z^*(z)| \simeq k.$$
By construction, we may pick the terms of the special
sequences close enough to $Y$, $Z$, $Y_*$, $Z_*$ respectively.
It is not difficult to check that we may
then perturb the  almost biorthonormal system of special sequences in such a way as to assume
that
$y \in Y$, $z \in Z$, $y^* \in Y_*$, $z^* \in Z_*$, and $y^*(z)=z^*(y)=0$, and preserving
the estimates on $\norm{y-z}\norm{y^*-z^*}$ and  $|y^*(y)-z^*(z)|$.

\

 When $X$ is reflexive,  the roles of
 $X$ and $X^*$ are interchangeable in the expression of $A$. 
Note that under reflexivity, the QHI property (\cite{F3} Corollary 4)  and the property of
having an unconditional basis are self-dual properties.

\

\subsection{FDD-block subspaces and FDD-block quotient of subspaces}
We shall prove that a
  dichotomy theorem  holds for quotient of subspaces which have a finite-dimensional decomposition
(or FDD) relative to a given Schauder basis (or even a FDD) of a given Banach space; they seem to be the
natural equivalent
of block-subspaces considered in Gowers' dichotomy.

\

An {\em interval of integers} is the intersection of $\N$ with a bounded interval of $\R$.
Two non empty intervals $E_1$ and $E_2$ are said to be {\em successive}, $E_1 < E_2$, when
$\max(E_1)<\min(E_2)$. A {\em successive partition} will be a sequence $(E_n)_{n \in \N}$
 of successive intervals
forming a partition of $\N$.

 Let $X$ be  a Banach space with a finite-dimensional decomposition denoted $(B_n)_{n \in \N}$. 
When $x=\sum_{n \in \N} b_n \in X$, with $b_n \in B_n$ for all $n \in \N$,
the {\em support of $x$} is the set 
$\{i \in \N: b_i \neq 0\}$. 
 The {\em range} of a vector is the smallest interval containing its support.
The support of a subspace $Y$  of $X$  is the smallest set 
containing the supports of all vectors
of $Y$. The range of $Y$  is the smallest interval containing the support of $Y$.
Two finitely supported subspaces $F$ and $G$ of $X$ with non-empty supports are {\em successive} when
$ran(F)<ran(G)$.

\

An {\em FDD-block subspace} of $X$ is
 an infinite sum $\sum_{n \in \N}F_n$ of finitely supported (possibly zero-dimensional) subspaces $F_n$
 of $X$,
such that $ran(F_n) \subset E_n, \forall n \in \N$, where $(E_n)_{n \in \N}$ is a successive partition.
Therefore an FDD-block subspace is finite-dimensional or equipped with the FDD $(F_n)_{n \in I}$, 
where $I=\{n: F_n \neq \{0\}\}$.

An {\em FDD-block quotient} of $X$ is the quotient of $X$ by some FDD-block subspace
$Y=\sum_{n \in \N}G_n$. An FDD-block quotient is finite-dimensional or
equipped with the FDD $(C_n)_{n \in I}$ corresponding to the successive partition $(E_n)_{n \in \N}$
associated to $Y$, that is,
where $C_n=([B_i,i \in E_n]+Y)/Y$ for all $n$, and
$I=\{n: C_n \neq \{0\}\}$.  Note that the space $X$ is an FDD-block quotient of itself.

An {\em FDD-block quotient of subspace} of $X$ is a quotient of subspace of $X$ of the form
$\sum_{n \in \N}F_n/\sum_{n \in \N}G_n$, where $G_n \subset F_n \subset [B_i, i \in E_n]$ for all $n$,
where $(E_n)_{n \in \N}$ is a successive partition.  The space
$\sum_{n \in \N}F_n/\sum_{n \in \N}G_n$ is naturally seen as an FDD-block subspace of $X/\sum_{n \in \N}G_n$,
when $X/\sum_{n \in \N}G_n$ is equipped with the FDD corresponding to $(E_n)_{n \in \N}$.
 
\

It is therefore clear that any FDD-block subspace (resp. quotient of subspace) of an FDD-block subspace (resp.
quotient of subspace) of $X$ is again
an FDD-block subspace (resp. quotient of subspace) of $X$. Note also that by classical results, 
any subspace of $X$ contains, for any $\epsilon>0$, an $1+\epsilon$-isomorphic
copy of a
block subspace, and therefore of an FDD-block subspace (however the similar result concerning
QS-spaces doesn't seem to be clear). Considering
FDD-block quotient of subspaces to study the structure
 of the class of QS-spaces is a natural counterpart
of considering block-subspaces to study the structure of the class of subspaces.

\begin{prop}\label{FDDQHI} Let $X$ be a Banach space with a finite-dimensional decomposition.
The following propositions are equivalent:

i) no FDD-block quotient of subspace of $X$ is decomposable,

ii) for any infinite codimensional 
FDD-block subspace $Y$ of $X$, the quotient $X/Y$ is hereditarily indecomposable,

iii) whenever $Y=\sum_{n \in \N}F_n/\sum_{n \in \N}G_n$
and $Y'=\sum_{n \in \N}F_n^{\prime}/\sum_{n \in \N}G_n$ are infinite dimensional FDD-block quotient of subspaces
of $X$ with a same successive partition, the sum
$Y+Y'$ is not direct in  $X/(\sum_{n \in \N}G_n)$.

\

When $X$ satisfies i) ii) iii) we shall say that
 $X$ is quotient hereditarily indecomposable restricted to FDD-block subspaces, or in short,
has the restricted QHI property.
\end{prop}

\pf ii) implies i) is immediate. If iii) is false then the FDD-block quotient of
subspace 
$Y+Y'=\sum_{n \in \N}(F_n+F_n^{\prime})/\sum_{n \in \N}G_n$ is decomposable, contradicting i). Finally, 
assume ii) is false,
i.e.
 $Z/W \oplus Z'/W$ forms a direct sum of infinite dimensional subspaces in $X/W$, for 
 some infinite codimensional
FDD-block subspace $W=\sum_{n \in \N}G_n$
 and some subspaces $Z$ and $Z'$, and let $(E_n)_{n \in \N}$ be a successive partition associated
to $W$.
We may up to a perturbation find sequences $(z_n)_{n \in \N}$ and $(z^{\prime}_n)_{n \in \N}$,
and a partition $(N_n)_{n \in \N}$ of $\N$ into successive intervals, 
such that for all $n \in \N$,
$ran(z_n,z^{\prime}_n) \subset \cup_{i \in N_n}E_i$, and such that $d(z_n,Z)$ and $d(z_n^{\prime},Z')$ converge to $0$ 
sufficiently fast so that $([z_n]_{n \in \N}+W)/W \oplus ([z^{\prime}_n]_{n \in \N}+W)/W$ is still direct
in $X/W$. Let for all $n \in \N$, $H_n=\sum_{i \in N_n}G_i$; we have therefore obtained
that $(\sum_{n \in \N} (H_n+[z_n]))/\sum_{n \in \N}H_n$ and
$(\sum_{n \in \N} (H_n+[z_n^{\prime}]))/\sum_{n \in \N}H_n$ form a direct sum, with successive partition
$(\cup_{i \in N_n}E_i)_{n \in \N}$, contradicting iii).
\pff

\

FDD-block quotient of subspaces still capture enough information about the structure of the space:
a space which has the restricted QHI property is in particular
 hereditarily indecomposable by ii), and by i), any of its infinite-dimensional FDD-block quotient of subspaces has again
the restricted QHI property.  The next proposition also shows that the restricted QHI property has similar 
self-dual properties as the QHI property.

\begin{prop}\label{QHIreflexif} Let $X$ be a Banach space with a shrinking 
finite-dimensional decomposition, such that $X^*$ has
the restricted
 QHI property. Then $X$ has the restricted QHI property.
\end{prop} 

\pf Let $Y=\sum_{n \in \N}F_n/\sum_{n \in \N}G_n$ be an infinite dimensional
 FDD-block quotient of subspace of $X$, with successive partition $(E_n)_{n \in \N}$.
 For each $n \in \N$, let $X_n^*$ be the 
space of vectors in $X^*$ with range included in $E_n$.
Then $$Y^*=(\sum_{n \in \N}G_n)^{\perp}/(\sum_{n \in \N}F_n)^{\perp}
=(\sum_{n \in \N}(G_n^{\perp}\cap X_n^*))/(\sum_{n \in \N}(F_n^{\perp}\cap X_n^*)),$$
since $(E_n)_{n \in \N}$ is a partition of $\N$.
So $Y^*$ is an FDD-block quotient of subspace of $X^*$. 
Therefore according to the first characterization in Proposition \ref{FDDQHI}, if $X$ does not
have the restricted QHI property, then  $X^*$ does not have the restricted QHI
property. \pff

\

In consequence, we note that if $X$ is a reflexive Banach space with the 
 restricted QHI property, then $X$ has HI dual and $X$ is saturated
with subspaces with HI dual. Indeed
every FDD-block subspace of $X$ has HI dual.  

\

We are now in position to state the result of this paper.

\begin{theo}\label{theo1}
Every Banach space has a quotient of subspace $Y$ with one of the two following properties, which
are mutually exclusive and both possible:

i) $Y$ has an unconditional basis,

ii) $Y$ has the restricted QHI property.
\end{theo}

We give a few comments on the reasons we needed to impose a restriction on the QHI property. Our proof
is based on some method of "combinatorial forcing", see Todorcevic's course \cite{A} about this. 
This will enable us to prove, up to some approximation, a general dichotomy result
 for closed properties of FDD-block quotient of 
subspaces, seen as sequences of finite dimensional successive QS-blocks
 (this will be defined precisely in the 
next section), with the product of the discrete topology
on the set of QS-blocks. This applies more or less directly to obtain Theorem \ref{theo1}.

 As we see them, these methods
 rely on defining 
 infinite sequences of elements which
may be correctly approximated by finite sequences; 
a notion of successivity
is needed, i.e. finite sequences are extended in infinite sequences
in a way that does not "affect" the properties implied by the finite part.
 
Our proof was inspired
by a simplification by B. Maurey of this method in the case of block-subspaces of a space with
 a Schauder basis, where a less restrictive 
setting may be used, based on replacing $X$ by a countable dense subset \cite{M2}.

 It
didn't seem possible to repeat exactly Maurey's proof to study QS spaces. Therefore we needed to restrict our study
to particular QS spaces which may be canonically associated to infinite sequences
 of "finite dimensional blocks"
 which are successive in some sense. For technical reasons, the countable dense subset must be replaced
 by a net whose intersection
with the set of predecessors of a given block is always finite. Up to perturbations, the restriction to a net 
is not essential, but the need for some notion of successivity seems to be, and this justifies that
we could not obtain "quotient hereditarily indecomposable" in the second part of the conclusion
of Theorem \ref{theo1}. Actually some examples indicate that FDD-block quotient of subspaces may behave differently
from general quotient of subspaces. We refer to the final section about this fact.

\section{Proof of the theorem}

To prove Theorem \ref{theo1}, we may consider a Banach space $X$ with a Schauder
basis $(e_n)$. We  denote by
$(e_n^*)$ its dual basis, and by $X_*$  the closed linear span of $(e_n^*)_{n \in \N}$.
 We may also assume, up to renorming, that the basis is
bimonotone. We shall consider supports and ranges of vectors, or subspaces, of $X$ and of $X_*$, with respect
to these canonical bases.

We choose to represent blocks forming
 quotient of subspaces of $X$ as pairs formed by
 a finite-dimensional subspace $F$  of $X$ and of 
a finite dimensional subspace of $F_*$ of  $X_*$, with $F_*^{\perp}
 \cap [e_n, n \in ran(F,F_*)] \subset F$. 
Pairs $(F,G)$ of finite-dimensional subspaces of $X$ with $G \subset F$ would also have been 
a possible representation. Our choice
will save us from some technicalities (successive pairs in our setting are pairs whose supports are
necessarily immediately successive). It will also preserve, in our proofs,
 the symmetry between the roles played by $X$ and $X^*$ in the 
reflexive case. This symmetry is apparent in our main result, and we felt it worth to be emphasized
in our demonstration.
 
\subsection{Blockings of QS-pairs}

If $Y \subset X$, and $Y_* \subset X_*$, the range of $(Y_*,Y)$ is the smallest interval
containing the ranges of $Y$ and of $Y_*$.
The set of finitely supported subspaces of $X$ is denoted
$F(X)$, of finitely supported subspaces of $X_*$ is denoted
$F(X_*)$.
A {\em QS-block (or block)}  is a pair $(F_*,F) \in F(X_*) \times F(X)$, therefore
$E:=ran(F_*,F)$ is finite,  such that $F_*^{\perp} \cap [e_n, n \in E] \subset F$.
The set of blocks is denoted ${\cal F}(X)$.
The {\em dimension} of $(F_*,F)$ is the dimension of $F/(F_*^{\perp} \cap [e_n, n \in E])$.
Two blocks $(F_*,F)$ and $(G_*,G)$ 
are said to be {\em successive} if $\min(ran(G_*,G))=\max(ran(F_*,F))+1$, and we write
$(F_*,F)<(G_*,G)$ (note the technical difference
with the usual notion of successivity).

We note that when ${\cal Y}=(Y_{n*},Y_n)_{n \in \N}$ is a sequence of successive 
blocks whose ranges partition $\N$ (equivalently, such that $\min(ran(Y_{1*},Y_1)=1$), the spaces
$Y=\sum_{n \in \N}Y_n$ and $Y_*=\sum_{n \in \N}Y_{n*}$ satisfy $Y_*^{\perp} \subset Y$. 
We shall then say that $(Y_*,Y)$ is the {\em QS-pair associated to ${\cal Y}$}, and that
${\cal Y}$ is {\em infinite dimensional} to mean that the QS-space $Y/Y_*^{\perp}$
is infinite dimensional. Note that the space  $Y/Y_*^{\perp}$ is a block-FDD quotient of subspace of $X$.

\

If $(F_*,F)$ is a block and
$(Y_{n*},Y_n)_{n \in I}$ is a finite or infinite sequence of successive blocks, 
and if there exists an interval $E \subset I$ such that
$F \subset \sum_{n \in E}Y_n$ and
$F_* \subset \sum_{n \in E}Y_{n*}$, then we shall say
that $(F_*,F)$ {\em is a block of} $(Y_{n*},Y_n)_{n \in I}$.

We now define a relations of "blocking" between sequences of successive blocks.

\begin{defi}
Let $(Y_{n*},Y_n)_{n}$ and $(Z_{i*},Z_i)_{i}$ be finite or infinite sequences of 
successive blocks. If for any $i$, $(Z_{i*},Z_i)$ is a block of $(Y_{n*},Y_n)_{n}$, then we 
shall say that $(Z_{i*},Z_i)_i$ is a blocking of $(Y_{n*},Y_n)_{n}$.

If ${\cal Z}=(Z_{i*},Z_i)_{i \in \N}$ and ${\cal Y}=(Y_{n*},Y_n)_{n \in \N}$ are infinite sequences of successive blocks
whose ranges partition $\N$, then we shall write
${\cal Z} \leq {\cal Y}$ to mean that
${\cal Z}$ is a blocking of ${\cal Y}$. This means 
that there exists a partition $\{N_i,i \in \N\}$ of $\N$ in successive intervals such
that, for all $i \in \N$, $ran(Z_{i*},Z_i)=\cup_{n \in N_i} ran(Y_{n*},Y_n)$ and
$(Z_{i*},Z_i)$ is a block of $(Y_{n*},Y_n)_{n \in N_i}$.

\end{defi}

We note that $\leq$ is an order relation. 
Clearly, when ${\cal Z} \leq {\cal Y}$, the associated 
QS-pairs $(Z_*,Z)$ and $(Y_*,Y)$
satisfy $(Z_*,Z) \subset (Y_*,Y)$.

For any two  sequences ${\cal Y}$ and ${\cal Z}$, we define 
$$A({\cal Y},{\cal Z})=A((Y_*,Y),(Z_*,Z)),$$
where $(Y_*,Y)$ and $(Z_*,Z)$ are the associated QS-pairs. 

\begin{lemm}\label{blabla}
Let ${\cal Y}=(Y_{n*},Y_n)_{n \in \N}$ be an infinite dimensional,
 successive sequence of blocks  whose ranges partition $\N$.
 Let $(Y_*,Y)$ be the associated QS-pair.
Assume that $A({\cal U},{\cal V})=0$ whenever ${\cal U},{\cal V} \leq {\cal Y}$ are infinite
dimensional, successive sequences of blocks whose  ranges are equal and partition $\N$.
Then $Y/Y_*^{\perp}$ has the restricted QHI property.
\end{lemm}

\pf The proof is based on the natural identification between sequences of blocks
of $Y/Y_*^{\perp}$ with its natural finite-dimensional decomposition,
 and sequences of blocks of $X$ which are blockings of $(Y_*,Y)$.  
Indeed consider
two infinite dimensional block-FDD quotient of subspaces of $Y/Y_*^{\perp}$
 which are of the form
 $Z=\sum_{n \in \N}F_n/\sum_{n \in \N}G_n$
and $Z'=\sum_{n \in \N}F_n^{\prime}/\sum_{n \in \N}G_n$, with successive partition $(E_n)_{n \in \N}$.
By definition  
for all $n \in \N$, $F_n \subset (\sum_{k \in E_n}Y_k+Y_*^{\perp})/Y_*^{\perp}$,
and let $I_n=ran(\sum_{k \in E_n}Y_k)$.
Therefore we may find $A_n,B_n$ such that
$$(\sum_{k \in E_n}Y_{k*})^{\perp} \cap [e_i, i \in I_n]
 \subset B_n \subset A_n \subset \sum_{k \in E_n}Y_n,$$
and such that $G_n=(B_n+Y_*^{\perp})/Y_*^{\perp}$ and 
$F_n=(A_n+Y_*^{\perp})/Y_*^{\perp}$. We define some subspaces $A_n^{\prime}$ associated to the spaces $F_n^{\prime}$ in a similar way.

We therefore have the identification 
$$Z = \overline{\sum_{n \in \N}A_n+Y_*^{\perp}}/\overline{\sum_{n \in \N}A_n+Y_*^{\perp}}=
\sum_{n \in \N}A_n/\sum_{n \in \N}B_n,$$ 
which is by construction a block FDD quotient of subspace of
$X$ corresponding to a blocking of ${\cal Y}$.
Indeed, let $B_{n*}=B_n^{\perp} \cap [e_i, i \in I_n]$, and let
${\cal Z}=(B_{n*},A_n)_{n \in \N}$, then the associated QS-space is
$\sum_{n \in \N}A_n/(\sum_{n \in \N}A_{n*})^{\perp}=\sum_{n \in \N}A_n/\sum_{n \in \N}B_n$.

We have the similar identification for $Z'$ and let
${\cal Z}'=(B_{n*},A_n^{\prime})_{n \in \N}$. Since  ${\cal Z} \leq {\cal Y}$ and
${\cal Z}' \leq {\cal Y}$, it follows that $A({\cal Z},{\cal Z}')=0$. This means
that the spaces $\sum_{n \in \N}A_n/\sum_{n \in \N}B_n$ and
$\sum_{n \in \N}A_n^{\prime}/\sum_{n \in \N}B_n$ do not form a direct sum in
$Y/\sum_{n \in \N}B_n$, and therefore
$Z$ and $Z'$ do not form a direct sum in the space $(Y/Y_*^{\perp})/(\sum_{n \in \N}G_n)$. Therefore iii) is satisfied
in Proposition \ref{FDDQHI}.\pff

\

Before stating more definitions, we need to
realize a reduction to a net ${\cal R}$ of blocks with some finiteness property which will be crucial for
our combinatorial method. 

For $F,G$ in $F(X)$, we let
$d_{H}(F,G)$  
be the Hausdorff distance between the unit spheres $S_F$ of $F$ and
$S_G$ of $G$, $d_H(F,G)=\max_{x \in S_F}d(x,S_G) \vee \max_{y \in S_G}d(y,S_F)$. 
Modifying a definition from \cite{F5}, we define a distance $d$ on $F(X)$ by 
$$d(F,G)=\min(1,2k\sqrt{k}d_H(F,G))$$ if $\dim F=\dim G=k$ and $ran(F)=ran(G)$, and
$d(F,G)=1$ otherwise.

We finally define a distance $\delta$ on ${\cal F}(X)$ by
$$\delta((F_*,F),(G_*,G))=\max(d(F,G),d(F_*^{\perp} \cap X_0,
G_*^{\perp} \cap X_0)),$$
when $ran(F_*,F)=ran(G_*,G)$ and $X_0=[e_i, i \in ran(F_*,F)]$, and we let
$\delta((F_*,F),(G_*,G))=1$ otherwise.

\

The critical result concerning this distance is contained in the next lemma.

\begin{lemm}\label{perturbations}
Let
$0<\epsilon<1$ and let
$(\delta_n)_n$ be a positive sequence such that $\sum_{n \in \N}\delta_n
\leq \epsilon$. Let
$(F_{n*},F_n)_{n
\in
\N}$ and
$(G_{n*},G_n)_{n
\in
\N}$ be successive sequences of blocks
such that for all $n \in \N$,
$\delta((F_{n*},F_n),(G_{n*},G_n)) \leq \delta_n$, and,
 let, for $n \in
\N$,
$X_n$ be the space $[e_i, i 
\in ran(F_{n*},F_n)]$.
Then there exists a map $T:\sum_{n \in \N}F_n \rightarrow \sum_{n \in
\N}G_n$ such that
$T(F_n)=G_n$ and $T(F_{n*}^{\perp}\cap X_n)=G_{n*}^{\perp} \cap X_n$
for all $n \in \N$, and such that for any $x \in \sum_{n \in \N}F_n$,
$\norm{Tx-x} \leq \epsilon\norm{x}$.
\end{lemm}

\pf Let $k=\dim F_1=\dim G_1$ and let $l=\dim F_{1*}^{\perp}\cap X_1
=\dim G_{1*}^{\perp}\cap X_1$. By classical results,
the Banach-Mazur distance of $F_1$ to $l_2^k$
is at most $\sqrt{k}$, so we may pick a normalized basis
$f_1,\ldots,f_k$ of $F_1$ such that $f_1,\ldots,f_l$ is a basis of
$F_{1*}^{\perp}\cap X_1$ and which has basis constant at most $\sqrt{k}$.
By the expression of $\delta$, we have that $d_H(F_{1*}^{\perp} \cap X_1,
G_{1*}^{\perp} \cap X_1) \leq \delta_1/2k\sqrt{k}$, therefore
for $1 \leq i \leq l$, there exists some $g_i \in G_{1*}^{\perp} \cap X_1$
with $\norm{g_i-f_i} \leq \delta_1/2k\sqrt{k}$.
Likewise we find for $l<i \leq k$ some $g_i \in G_1$ with
the same condition on $\norm{g_i-f_i}$.

By \cite{LT} Prop. 1.a.9, $(g_i)_{1 \leq i \leq k}$ is a basis of
$G_1$, and furthermore, if $T_1:F_1 \rightarrow G_1$ is defined
by $T_1(f_i)=g_i$ for all $1 \leq i \leq k$, we have, for any $x \in F_1$,
$x=\sum_{i=1}^k a_i f_i$,
$$\norm{T_1x-x} \leq \sum_{i=1}^k |a_i|\norm{f_i-g_i} 
\leq 2\sqrt{k}\norm{x} k(\delta_1/2k\sqrt{k}) \leq \delta_1 \norm{x}.$$

Repeating this construction on each $F_n$,  let $T_n$ be the associated
map from $F_n$ onto $G_n$ with $T_n(F_{n*}^{\perp} \cap X_n)=
G_{n*}^{\perp} \cap X_n$, and let $T$ be defined on $\sum_{n \in \N}F_n$ 
by $T_{|F_n}=T_n$ for all $n \in \N$.
We have
for any $x=\sum_{n \in \N}x_n, x_n \in F_n$,
$$\norm{Tx-x} \leq \sum_{n \in \N}\norm{T_n x_n-x_n} \leq
\sum_{n \in \N}\delta_n \norm{x_n} \leq \epsilon \norm{x},$$
by bimonotonicity of the basis. \pff

\

For $N \in \N$, we let 
${\cal F}_N(X)$ be the set of elements $(F_*,F)$
 of ${\cal F}(X)$ such that $\max(ran(F_*,F))=N$.
Fixing $(\delta_n)_{n \in \N}$  a decreasing positive sequence such that $\delta_n \leq 2^{-n}$ for every
 $n \in \N$,
we  define 
${\cal R} \subset {\cal F}(X)$ satisfying the following
properties:

\

i) 
${\cal R} \cap {\cal F}_N(X)$ is a finite $\delta_N$-net for ${\cal F}_N(X)$,

ii) whenever $(F_{1*},F_1)<\cdots<(F_{k*},F_k)$ belong to ${\cal R}$, it follows that
$(F_{1*}+\ldots+F_{k*},F_1+\ldots+F_k)$ belongs to ${\cal R}$.

iii) for any $(F_*,F) \in {\cal R} \cap {\cal F}_N(X)$,
${\cal R} \cap {\cal F}_{F_*,F}$ is a $\delta_N$-net
for ${\cal F}_{F_*,F}$, where
 ${\cal F}_{F_*,F}$ denotes $\{(G_*,G) \in {\cal F}_N(X): (G \subset F) \wedge (G_* \subset F_*)\}$.

iv) for any $(F_*,F) \in {\cal R} \cap {\cal F}_N(X)$,
${\cal R} \cap {\cal F}_F^{F_*}$ is a $\delta_N$-net
for ${\cal F}_{F}^{F_*}$, where
${\cal F}_F^{F_*}:=\{(G_*,F) \in {\cal F}_N(X): G_* \subset F_*\}$,

v) for any $(F_*,F) \in {\cal R} \cap {\cal F}_N(X)$,
${\cal R} \cap {\cal F}_{F_*}^F$ is a $\delta_N$-net
for ${\cal F}_{F_*}^F$, where
 ${\cal F}_{F_*}^F:=\{(F_*,G) \in {\cal F}_N(X): G \subset F\}$,

vi) if $(F_*,F) \in {\cal R}$ then 
$(F^{\perp}\cap [e_i^*,i \in E],F_*^{\perp} \cap [e_i,i \in E]) \in {\cal R}$,
where $E=ran(F_*,F)$. 

\

An {\em ${\cal R}$-block} will denote a block in ${\cal R}$.
In the following, blocks will always be  ${\cal R}$-blocks, unless specified otherwise.

\

We denote by $QS^{<\omega}(X)$  (resp. $QS_0^{<\omega}(X)$)
the set of finite sequences of successive $\cal R$-blocks 
$(F_{n*},F_n)_n$ (resp. for which $\min(ran(F_{1*},F_1))=1$).

The set $QS^{\omega}(X)$ (resp. $QS_0^{\omega}(X)$)
is the space of infinite sequences of successive $\cal R$-blocks
${\cal Y}=(Y_{n*},Y_n)_n$ (resp. for which $\min(ran(Y_{1*},Y_1))=1$).
If $(Y_{n*},Y_n)_n$ is  an element of $QS_0^{\omega}(X)$, the {\em partition} of
$(Y_{n*},Y_n)_n$ is the sequence $(ran(Y_{n*},Y_n))_{n \in \N}$, which forms a
partition of $\N$.
 The space $\sum_{n \in \N}Y_n$ will be denoted $Y$,
and $Y_*$ will denote $\sum_{n \in \N}Y_{n*}$. As was already observed, the relation
 $Y_*^{\perp} \subset Y$
ensures that $Y/Y_*^{\perp}$ is a block-FDD quotient of subspace of $X$.
We let $QS(X) \subset QS_0^{\omega}(X)$ be the set of sequences which are infinite dimensional, that is
such that
the QS-space $Y/Y_*^{\perp}$ is infinite dimensional.

\

If $E$ is an interval of integers, and $(Y_{n*},Y_n)_{n \in I}$ is a finite or infinite sequence
 of successive $\cal R$-blocks, we shall say that $(Y_{n*},Y_n)_{n \in I}$ {\em is well-placed with respect to
$E$} if there exists $m \in I$ such that $\min(ran(Y_{m*},Y_m))=\max E+1$.
The set of sequences of $QS(X)$ which are well-placed with respect to $E$ is denoted
$QS_E(X)$.

\

We now define a relation of "tail blocking" on $QS(X)$.

\begin{defi}
Let ${\cal Z},{\cal Y} \in QS(X)$. If $E$ an interval of $\N$,  and $(Z_{i*},Z_i)_{i \geq p}$ is 
 a blocking of $(Y_{n*},Y_n)_{n \geq m}$, with
$\min(ran(Z_{p*},Z_p))=\min(ran(Y_{m*},Y_m))=\max E+1$, then we shall write
that ${\cal Z} \leq^E {\cal Y}$. 
\end{defi}

\

Note that if ${\cal Z} \leq^E {\cal Y}$ then it follows necessarily that
${\cal Z}$ and ${\cal Y}$ are well-placed with respect to $E$.
 It is also clear that $\leq^E$ a preorder relation, and that ${\cal W} \leq^E {\cal Y}$ whenever
${\cal W}$ and ${\cal Y}$ are well-placed with respect to $E$ and
${\cal W}  \leq {\cal Y}$.

We shall need the following easy lemma.

\begin{lemm}\label{easylemma} Let $E$ be an interval of $\N$, ${\cal Y}, {\cal Z} \in QS_E(X)$. 
Assume ${\cal Z} \leq^E {\cal Y}$.
Then there exists  ${\cal W} \in QS_E(X)$ such that ${\cal W} \leq {\cal Y}$
and ${\cal W} \leq^E {\cal Z}$.
\end{lemm}

\pf Let
$\min(ran((Y_{m*},Y_m)))=\max E+1=\min(ran((Z_{p*},Z_p)))$ for some $m,p$.
We define $(W_{n*},W_n)=(Y_{n*},Y_n)$ if $n<m$ and
$(W_{n*},W_n)=(Z_{(n-p+m)*},Z_{n-p+m})$ if $n \geq m$.\pff

\begin{defi} Let $P \subset QS_E(X)$. We say that $P$ is {$\leq^E$-hereditary} if
 whenever ${\cal Y} \in P$ and ${\cal Z} \leq^E {\cal Y}$, then
${\cal Z} \in P$. We say that $P$ is
$\leq^E$-large if it is $\leq^E$-hereditary  and
 whenever ${\cal Y} \in QS_E(X)$, there exists ${\cal Z} \leq {\cal Y}$ such that
${\cal Z} \in P$.
\end{defi}

\subsection{A game for QS-pairs}

Our proof will be based on an "oriented QS-pairs" Gowers game
 $G_{\cal A}^{\cal Y}$ associated to some 
subset $\cal A$ of $QS(X)\times \{-1,1\}^{\omega}$ and to some ${\cal Y} \in QS(X)$, and
 defined as follows.
Player 1 plays some  ${\cal W}_1 \leq {\cal W}$.
Player 2 plays some sign $\epsilon_1 \in \{-1,1\}$, and some block
$(U_{1*},U_1)$ which is a block of ${\cal W}_1$ with $\min(ran(U_{1*},U_1))=1$.

At step $n$, Player 1 plays some ${\cal W}_n \leq {\cal Y}$ which is well-placed with
respect to $ran(U_{n-1*},U_{n-1})$. Player 2 plays some sign $\epsilon_n \in \{-1,1\}$, and some
 block $(U_{n*},U_n)$ of ${\cal W}_n$ which is successive with respect to
  $(U_{n-1*},U_{n-1})$.

Player 2 wins the game if he produced an infinite sequence $(U_{n*},U_n,\epsilon_n)_n$ which is in ${\cal A}$.

\

In our application we shall use this game for the set ${\cal A}_{\delta}$ associated to some $\delta>0$ and 
defined as the set of $(U_{n*},U_n,\epsilon_n)_n$ such that
there exists $n \in \N$, there exists $u_k \in U_k$, $u_k^* \in U_{k*}$,
 $1 \leq k \leq n$, such 
that $$\norm{\sum_{k=1}^n u_k}\norm{\sum_{k=1}^n u_k^*}
<\delta |\sum_{k=1}^n \epsilon_{k-1} u_k^*(u_k)|,$$
where $\epsilon_0=1$ is fixed.

\

A {\em state} $s$ will be an element of $QS^{<\omega}(X) \times\{-1,1\}^{<\omega}$, where 
the two sequences are
of equal length denoted $|s|$. The set of states will be denoted $S$.
When $\cal Y$ is well-placed with respect to $(U_{i*},U_i)_{i<k}$ and $(\epsilon_i)_{i<k}$ is
a sequence of signs, we define in an obvious way the game
$G_{{\cal A}}^{\cal Y}(s)$, where
$s$ is the state $(U_{n*},U_n,\epsilon_n)_{n < k}$: just rename the steps $1,2,\ldots$ in the new game  step
$k,k+1,\ldots$ and then apply the same definition as above; this is the game
$G_{\cal A}^{\cal Y}$ starting from position $s$. 

If $s=(U_{n*},U_n,\epsilon_n)_{n < k}$, then
$ran(s)$ will denote $ran((U_{n*},U_n)_{n < k})$, and to simplify the notation we also
 let $QS_s(X)$ stand for
$QS_{ran(s)}(X)$, $\leq^s$ stand for $\leq^{ran(s)}$, "successive to $s$" mean "successive
to $(U_{k-1*},U_{k-1})$". 

\

In the following, we fix some subset ${\cal A}$ of $QS(X)\times \{-1,1\}^{\omega}$.
Our next definition is the first step of  the method of "combinatorial forcing" on $QS(X)$.

\begin{defi} Let $s$ be a state, and let ${\cal Y} \in QS_s(X)$.

 The state $s$ accepts ${\cal Y}$ if Player 2 has a winning strategy
for the game $G_{{\cal A}}^{\cal Y}(s)$.

 The state $s$ rejects ${\cal Y}$ if it
accepts no  ${\cal Z} \leq {\cal Y}$.

The state $s$ decides ${\cal Y}$ if it accepts or
 rejects ${\cal Y}$.
\end{defi}

\begin{lemm}\label{decides}
Let $s$ be a state. 

- the set of ${\cal Y}$ in $QS_s(X)$
such that $s$ accepts ${\cal Y}$ 
(resp. rejects ${\cal Y}$) is $\leq^s$-hereditary.

- the set of ${\cal Y}$ in $QS_s(X)$
such that $s$ decides ${\cal Y}$ 
 is $\leq^s$-large.
\end{lemm}

\pf 
Assume $s$ accepts ${\cal Y}$.
Let ${\cal Z}$ be such that
${\cal Z} \leq_s {\cal Y}$.  
Let at step $n$, ${\cal W}={\cal W}_n \leq {\cal Z}$ be a move for Player 1. 
By Lemma \ref{easylemma}, we may find
 ${\cal V} \leq {\cal Y}$ with ${\cal V} \leq^s {\cal W}$, in particular
 ${\cal V}$ is well-placed
with respect to $ran(s)$. Therefore ${\cal V}_n={\cal V}$ is an admissible move 
for Player 1.  Since $s$ accepts, a move $(U_{n*},U_n,\epsilon_n)$ for
Player 2 is prescribed by the winning strategy for $G_{{\cal A}}^{\cal Y}(s)$.
This move is admissible for Player 2 in $G_{{\cal A}}^{\cal Z}(s)$, since
$(U_{n*},U_n)$ is successive to $s$ and therefore is a block of
  ${\cal W}$.
We have therefore described a winning strategy for Player 2 in the game
$G_{{\cal A}}^{\cal Z}(s)$, which means that $s$ accepts ${\cal Z}$.

Assume now  that $s$ rejects ${\cal Y} \in QS_s(X)$ while
it does not reject   ${\cal Z} \leq^s {\cal Y}$. We may assume that 
$s$ accepts ${\cal Z}$. We get a contradiction by using Lemma \ref{easylemma} to find
some element ${\cal W} \in QS_s(X)$ such that ${\cal W} \leq {\cal Y}$
and ${\cal W} \leq^s {\cal Z}$.

It follows from this that the set of ${\cal Y}$ in $QS_s(X)$
such that $s$ decides ${\cal Y}$ 
 is $\leq^s$-hereditary. Finally if ${\cal Y} \in QS_s(X)$, either $s$ rejects
${\cal Y}$, or $s$ accepts some ${\cal Z} \leq {\cal Y}$; this implies
$\leq^s$-largeness.\pff

 \begin{lemm}\label{stabilizes}(stabilization principle) For any ${\cal W} \in QS(X)$,
there exists ${\cal Y} \leq {\cal W}$ such that
whenever $(Z_{n*},Z_n)_{n \leq k} \in QS_0^{<\omega}(X)$ is a blocking of ${\cal Y}$,
and $(\epsilon_n)_{n \leq k}$ is a sequence of signs, it follows that the state
$s=(Z_{n*},Z_n,\epsilon_n)_{n \leq k}$ decides ${\cal Y}$. 

Such a
${\cal Y}$ will be called stabilizing, and states associated to
blockings of ${\cal Y}$ will be said to be states blocking ${\cal Y}$.
\end{lemm}

\pf Let ${\cal W}$ be fixed in $QS(X)$. Let $n_1$ be such that
 $\dim(W_{n_1 *},W_{n_1}) \geq 1$. We let ${\cal Y}^1={\cal W}$ and let
$(Y_{1*},Y_1)=(\sum_{n \leq n_1}Y^1_{n*},\sum_{n \leq n_1} Y^1_n)$.
Assume $(Y_{k*},Y_k)_{k<n}$ and some ${\cal Y}^{n-1}$ in $QS_E(X)$ were constructed with
 $E=ran((Y_{n-1 *},Y_{n-1})$.
By the finiteness property of $\cal R$, and the $\leq^E$-largeness property of Lemma \ref{decides},
we may find some ${\cal Y}^n \leq {\cal Y}^{n-1}$, with ${\cal Y}^n \in QS_E(X)$, such that
for any finite sequence $(Z_{i*},Z_i)_{i \leq m}$ which is a blocking of
$(Y_{k*},Y_k)_{k<n}$ with $\max(ran(Z_{m*},Z_m))=\max E$, and
for any sequence of signs $(\epsilon_i)_{i \leq m}$, the state
$s=(Z_{i*},Z_i,\epsilon_i)_{i \leq m}$ decides ${\cal Y}^n$.
Let $m_n$ be such that
$\max(ran(Y^n_{m_n*},Y^n_{m_n})=\max E$ and $p_n$
be such that the associated subsequence $(Y^n_{i*},Y^n_i)_{m_n <i \leq p_n}$
 contains a term of dimension at least $1$.
Let $(Y_{n*},Y_n)$ be  
$(\sum_{m_n<i \leq p_n}Y^{n}_{i*},\sum_{m_n<i \leq p_n}Y^{n}_{i})$.

Repeating this by induction we have constructed an element of $QS(X)$ 
which satisfies the required property.
Indeed for any state $s$ blocking $\cal Y$, let
$n$ be such that $\max(ran(s))=\max(ran(Y_{n-1 *},Y_{n-1}))$. Then
$s$ decides ${\cal Y}^n$ and ${\cal Y} \leq^s {\cal Y}^n$, therefore
$s$ decides ${\cal Y}$.\pff

\

We now fix some stabilizing  ${\cal X}$ in $QS(X)$.
 Note that by Lemma \ref{decides}  and Lemma \ref{stabilizes},
whenever $s$ is a state blocking ${\cal X}$ and ${\cal Y} \leq_s {\cal X}$, we have that
$s$ accepts (resp. rejects) ${\cal X}$ if and only if
it accepts (resp. rejects) ${\cal Y}$.
In the following, we shall write $s$ accepts (resp. rejects), to mean
 that $s$ accepts
(resp. rejects) ${\cal X}$.

\

\begin{lemm}\label{sublemma}
Let $s \in S$ be a state blocking ${\cal X}$. If $s$ rejects, then
for any  ${\cal Y} \leq {\cal X}$ in $QS_s(X)$
 there
exists ${\cal Z} \leq {\cal Y}$ in $QS_s(X)$ such that
for any $(F_*,F)$ block of ${\cal Z}$ which is successive to $s$, and
any sign $\epsilon$, the state
$s^{\frown}(F_*,F,\epsilon)$ rejects.
 \end{lemm}

\pf 
Assume the conclusion is false. Let $n=|s|$.
There exists ${\cal Y} \leq {\cal X}$ in $QS_s(X)$, such that
for any ${\cal Z} \leq {\cal Y}$ in $QS_s(X)$,
 there is a block $(F_{n+1 *},F_{n+1})$ of ${\cal Z}$
successive to $s$ and  
  $\epsilon_{n+1} \in \{-1,1\}$
such that the state $s'=s^{\frown}(F^*_{n+1},F_{n+1},\epsilon_{n+1})$ accepts, and therefore
accepts $\cal Y$,
 that is Player 2 has a
winning strategy for $G_{\cal A}^{{\cal Y}}(s')$. Note that
$s'$ is a state blocking ${\cal X}$.
What we wrote means that Player 2 has a
winning strategy for $G_{\cal A}^{\cal Y}(s)$, in other words $s$ accepts $\cal Y$, that is $s$ accepts.
This is a contradiction. \pff

\

In the following $\emptyset$ denote the empty state.

\begin{lemm}\label{rejects}
Assume $\emptyset$ rejects. Then there exists ${\cal Y} \leq {\cal X}$
 such that any state blocking ${\cal Y}$  rejects.
\end{lemm}

\pf Let ${\cal Y}^0={\cal X}$. We build by induction
a sequence ${\cal Y}=(Y_{n*},Y_n)_{n \in \N}$ and a $\leq$-decreasing sequence $({\cal Y}^n)_{n \in \N}$
with ${\cal Y}^n \in QS_{E_n}(X)$, if $E_n=ran(Y_{i*},Y_i)_{i<n}$, and with $(Y_{n*},Y_n)$  a block of
${\cal Y}^n$ for each $n \in \N$, as follows.
Assume $(Y_{i*},Y_i)_{i<n}$  and 
$({\cal Y}^{i})_{i<n}$  were defined.
There are finitely many states $s$ with $\max(ran(s))=\max(E)$.
  Therefore applying
Lemma \ref{sublemma} a finite number of times, we obtain some
${\cal Y}^n \leq {\cal Y}^{n-1}$ in $QS_E(X)$ such that for any state $s$ with $\max(ran(s))=\max(E)$,
for any $(F_*,F)$ block of ${\cal Z}$ which is successive to $E$, and for
any sign $\epsilon$, the state
$s^{\frown}(F_*,F,\epsilon)$ rejects.
We define $(Y_{n*},Y_n)$ to be such a block $(F_*,F)$ of dimension at least $1$.

Whenever ${\cal U}=(U_{n*},U_n)_{n \in \N} \leq {\cal Y}$, we may easily check by induction
that for any sequence of signs $(\epsilon_i)_{i \leq n}$,
the state $(U_{i*},U_i,\epsilon_i)_{i \leq n}$ rejects. \pff

\

\begin{prop} Let $\cal A$ be a subset of $QS(X)\times \{-1,1\}^{\omega}$
 which is open as a subset of 
$({\cal F}(X) \times \{-1,1\})^{\omega}$ with the product of the discrete topology on ${\cal F}(X) \times
\{-1,1\}$.
 If for every
${\cal Y} \in QS(X)$, there exists ${\cal Z} \leq {\cal Y}$ and a sequence of signs $e$ 
such that $({\cal Z},e) \in {\cal A}$, then
there exists ${\cal Y} \in QS(X)$ such that Player 2 has a winning strategy in the game 
$G_{\cal A}^{\cal Y}$. \end{prop}

\pf If $\emptyset$ accepts then by definition,
 Player 2 has a winning strategy in the game 
$G_{\cal A}^{\cal Y}$ for some ${\cal Y}$. If $\emptyset$ rejects then, by 
Lemma \ref{rejects}, there exists ${\cal Y}$ of which any blocking state 
rejects, which implies that any  state blocking  ${\cal Y}$ is extendable as a sequence which is not in ${\cal A}$.
Since ${\cal A}$ is open, this means that no infinite sequence of successive blocks of ${\cal Y}$ and of signs
belongs 
to ${\cal A}$. \pff

\

Recall that for any $\delta>0$, we define ${\cal A}_{\delta}$  to be the set of
 $(U_{n*},U_n,\epsilon_n)_n$ such that
there exists $n \in \N$, and $u_k \in U_k$, $u_k^* \in U_{k*}$,
 $1 \leq k \leq n$, such 
that $$\norm{\sum_{k=1}^n u_k}\norm{\sum_{k=1}^n u_k^*}
<\delta |\sum_{k=1}^n \epsilon_{k-1} u_k^*(u_k)|,$$
where we put $\epsilon_0=1$.
This is an open subset of $({\cal F}(X) \times \{-1,1\})^{\omega}$.

\subsection{A dichotomy theorem on QS(X)}

If ${\cal Y} \in QS(X)$, with
$\dim(Y_{n*},Y_n)=1$ for all $n \in \N$, then we shall write
that ${\cal Y} \in QS_1(X)$. If ${\cal Y} \in QS_1(X)$, and for each $n \in \N$,
$\tilde{e}_n \in Y/Y_*^{\perp}$ is the class of some $e_n \in Y_n$ which is not in $Y_{n*}^{\perp}$,
  then we shall say that
$(\tilde{e}_n)$ is a {\em successive Schauder basis} of $Y/Y_*^{\perp}$. Note that all
successive Schauder bases of $Y/Y_*^{\perp}$ may be deduced from each other by homotheties on the span of
each
of their basic vectors.

\

In the next proposition, fixing $\delta>0$, we let ${\cal X}_{\delta}$ be a stabilizing subspace corresponding to ${\cal A}_{\delta}$,
and we write $s$ $\delta$-accepts (resp. $\delta$-rejects) to mean that
$s$ accepts (resp. rejects) ${\cal X}_{\delta}$ with respect to the set ${\cal A}_{\delta}$.

\begin{prop}\label{prop29}
If $\emptyset$ $\delta$-rejects, then there exists ${\cal Y} \in QS_1(X)$ with
 ${\cal Y} \leq {\cal X}_{\delta}$, such that
any successive basis of $Y/Y_*^{\perp}$  is unconditional  with constant $\delta^{-1}$. 
If $\emptyset$ $\delta$-accepts, then whenever ${\cal U}, {\cal V} \leq {\cal X}_{\delta}$ 
have identical partitions, $A({\cal U}, {\cal V})<\delta$.
\end{prop}

\pf If $\emptyset$ $\delta$-rejects, then consider ${\cal Y}=(Y_{i*},Y_i)_{i \in \N}$ 
given by Lemma \ref{rejects}, and 
write $E_i=ran(Y_{i*},Y_i)$. Without loss of 
generality we may assume that ${\cal Y}$ belongs to $QS_1(X)$.
Pick in each $Y_i$ some normalized $f_i$ such that $d(f_i,Y_{i*}^{\perp} \cap [e_n, n \in E_i])=1$. Fix $n$ and some signs
$(\epsilon_i)_{i \leq n}$, and recall that $\epsilon_0$=1. By the proof of Lemma \ref{rejects}, and since
${\cal A}_{\delta}$ is open, we have  that
for any  $(y_i^*,y_i) \in Y_{i*}\times Y_i$, $i \leq n$,

$$\norm{\sum_{k=1}^n y_k}\norm{\sum_{k=1}^n y_k^*}
\geq \delta |\sum_{k=1}^n \epsilon_{k-1} y_k^*(y_k)|.$$

Equivalently,  whenever $\norm{\sum_{k=1}^n y_k^*}=1$,
$$(\sum_{k=1}^n y_k^*)(\sum_{k=1}^n\epsilon_{k-1} y_k) \leq 
\delta^{-1}\norm{\sum_{k=1}^n y_k},$$
and therefore
$$\norm{\sum_{k=1}^n\epsilon_{k-1} y_k}_{\sum_{k \leq n}Y_{k*}} \leq
\delta^{-1}\norm{\sum_{k=1}^n y_k}.$$
Taking $y_k=\lambda_k f_k+z_k$, where $\lambda_k$ is a real number and $z_k$ is arbitrary in
$Y_{k*}^{\perp} \cap [e_i,i \in E_k]$, we obtain
$$\norm{\sum_{k=1}^n\epsilon_{k-1} \lambda_k f_k}_{\sum_{k \leq n} Y_{k*}} \leq
\delta^{-1}\norm{\sum_{k=1}^n \lambda_k f_k+z},$$
where $z \in (\sum_{k \leq n} (Y_{k*}^{\perp} \cap [e_i, i \in E_k])
=(\sum_{k \leq n}Y_{k*})^{\perp} \cap [e_i, i \in \cup_{k \leq n} E_k]$ is arbitrary.
 By duality in $[e_i, i \in \cup_{k \leq n} E_k]$, we conclude that
$$\norm{\sum_{k=1}^n\epsilon_{k-1} \lambda_k f_k}_{\sum_{k \leq n} Y_{k*}} \leq
\delta^{-1}\norm{\sum_{k=1}^n \lambda_k f_k}_{\sum_{k=1}^n Y_{k*}}.$$
Since $(\epsilon_i)_{1 \leq i \leq n-1}$ was arbitrary, we deduce
that  $(\tilde{e}_k)_{k \leq n}$ is $\delta^{-1}$-unconditional in 
$\sum_{k \leq n}Y_k/((\sum_{k \leq n}Y_{k*})^{\perp} \cap [e_i, i \in \cup_{k \leq n}E_k])$ for 
each $n$, and therefore in $Y/Y_*^{\perp}$ by bimonotonicity.

\

Assume $\emptyset$ $\delta$-accepts. Pick ${\cal U}, {\cal V} \leq {\cal X}_{\delta}$ which
have identical partitions. This will ensure that playing
${\cal U}$ or ${\cal V}$
is always an admissible move for Player 1. We therefore may
define a strategy
for Player 1 as follows. The first move is $\cal U$. Assuming Player $2$ picked
some $(Y_{k-1}^*,Y_{k-1},\epsilon_{k-1})$ at step $k-1$, Player $1$'s $k$-th move will be
$\cal U$ if $\epsilon_{k-1}=1$ and $\cal V$ if $\epsilon_{k-1}=-1$.
Opposing a winning strategy for Player 2, we therefore obtain
some $n \in \N$, and some sequences $(u_{i}^*,u_i)_{i \leq n}$ of pairs of vectors
and functionals, and $(\epsilon_i)_{i \leq n}$ of signs such that
$u_i \in U, u_{i}^* \in U_*$ if $\epsilon_{i-1}=1$ and
$u_i \in V, u_{i}^* \in V_*$ if $\epsilon_{i-1}=-1$, and with
$$\norm{\sum_{k=1}^n u_k}\norm{\sum_{k=1}^n u_k^*}
<\delta |\sum_{k=1}^n \epsilon_{k-1} u_k^*(u_k)|.$$
We let $u=\sum_{\epsilon_{i-1}=1}u_i \in U$,
$u^*=\sum_{\epsilon_{i-1}=1}u^*_i \in U_*$,
$v=-\sum_{\epsilon_{i-1}=-1}u_i \in V$,
$v^*=-\sum_{\epsilon_{i-1}=-1}u_i^* \in V_*$, and observe that $u^*(v)=v^*(u)=0$ and
$$\norm{u-v}\norm{u^*-v^*}
<\delta |u^*(u)-v^*(v)|.$$ \pff

\

\begin{theo}\label{tadam} Let $X$ be a Banach space with a Schauder basis. Then there exists
 a quotient of subspace $Y/Y_*^{\perp}$ of $X$, associated to some
${\cal Y}$ in $QS_1(X)$, which satisfies one of the two following properties, which are both possible and
mutually exclusive: 

i) $Y/Y_*^{\perp}$ has an unconditional basis,

ii) $Y/Y_*^{\perp}$ has the restricted QHI property.  
\end{theo}

\pf Fix as before a positive sequence $(\delta_n)_{n \in \N}$ with
$\delta_n \leq 2^{-n}$ for all $n$, and build by Lemma \ref{stabilizes} a $\leq$-decreasing sequence 
${\cal X}_n$ such that ${\cal X}_n$ is $\delta_n$-stabilizing for 
each $n$.
If, with the notation defined as the beginning of this subsection,
$\emptyset$ $\delta_n$-rejects ${\cal X}_n$ for some $n$, then we are done by Proposition \ref{prop29}.

Assume therefore that $\emptyset$ $\delta_n$-accepts ${\cal X}_n$ for all $n \in \N$.
 Let ${\cal Y} \in QS(X)$ be diagonal for the ${\cal X}_n$'s, i.e. such that
for any state $s$ blocking ${\cal Y}$, with
 $\max(ran(s))=\max(ran(Y_{n*},Y_n))$, we have that
${\cal Y} \leq^s {\cal X}_n$. This is easily constructed by induction.
We shall prove that $A({\cal U},{\cal V})=0$ for any ${\cal U},{\cal V} \leq {\cal Y}$ which
are sequences of successive blocks (not necessarily in ${\cal R}$) with equal ranges
forming  a partition
of $\N$. By Lemma \ref{blabla}, this will be enough to prove our result.

\

Fix $\epsilon>0$, and arbitrary ${\cal U},{\cal V} \in QS(X)$ (therefore formed of ${\cal R}$-blocks), with
 ${\cal U},{\cal V} \leq {\cal Y}$ and with identical partition denoted
$(E_n)_{n \in \N}$.  Let $m$ be large enough so that if
$p=\max(E_m)$ then $\delta_p<\epsilon$.
Denote by ${\cal X}=(X_{i*},X_i)_{i \in \N}$ the corresponding ${\cal X}_p$, 
by $F_i$ the range of $(X_{i*},X_i)$ and let $q$ be such that
$p=\max(F_q)$.

We let for $i \leq q$, $U^{\prime}_{i*}=X_{i*}$, and $U^{\prime}_i=
X_{i*}^{\perp} \cap [e_n, n \in F_i]$. For 
$i > q$ we let $(U^{\prime}_{i*},U^{\prime}_i)=(U_{m-q+i *},U_{m-q+i})$.
 We have therefore 
constructed an element ${\cal U}^{\prime}$ of $QS(X)$ which 
satisfies ${\cal U}^{\prime} \leq {\cal X}$ and ${\cal U}^{\prime} \leq^E {\cal U}$ for $E=[1,p]$.
We construct in the same way some
${\cal V}^{\prime} \leq {\cal X}$, ${\cal V}^{\prime} \leq^E {\cal V}$.

By Proposition \ref{prop29}, we may find
 $x,x^* \in {\cal U}^{\prime}$ and
$y,y^* \in {\cal V}^{\prime}$, with disjoint supports, with
$\norm{x-y}\norm{x^*-y^*}<\delta_p
|(x^*-y^*)(x+y)|$.
Let $P$ be the projection onto $[e_n, n>p]$, and
 $P_*$ be the projection onto $[e_n^*, n>p]$.
Note that $P(U') \subset U$ and $P_*(U^{\prime}_*) \subset U_*$, and the similar inclusions
hold for $V$ and $V_*$.
Let $u=Px, u^*=P_*x^*, v=Py, v^*=P_*y^*$.

By bimonotonicity of the basis, we observe that
$\norm{u-v} \leq \norm{x-y}$ and
$\norm{u^*-v^*} \leq \norm{x^*-y^*}$. 
On the other hand, writing $x=u+a, x^*=u^*+a^*, y=v+b, y^*=v^*+b^*$, we note that
$a \in (\sum_{i \leq q}X_{i*}^{\perp}$ while $a^* \in  \sum_{i \leq q}X_{i*}$, therefore
$a^*(a)=0$. Likewise, $b^*(b)=a^*(b)=b^*(a)=0$, and by disjointness of the ranges,
$u^*(a)=u^*(b)=v^*(a)=v^*(b)=a^*(u)=a^*(v)=b^*(u)=b^*(v)=0$.
Therefore
$$(x^*-y^*)(x+y)=(u^*-v^*)(u+v),$$
and we deduce that
$$\norm{u-v}\norm{u^*-v^*}<\delta_p
|(u^*-v^*)(u+v)|.$$
We have therefore proved that $a({\cal U},{\cal V}) <\epsilon$.

\

It remains to show that we may obtain the same results for general ${\cal U},{\cal V} \leq {\cal Y}$, i.e. 
successive sequences of  blocks which are not necessarily in ${\cal R}$. Fix $0<\epsilon<1/3$
 and let ${\cal U}$, ${\cal V}$ have the
same partition $(E_n)_{n \in \N}$. Let ${\cal U}^{\prime}$, ${\cal V}^{\prime}$ be sequences with
blocks in ${\cal R}$, such that $\delta((U_{n*},U_n),(U^{\prime}_{n*},U^{\prime}_n))<\delta_n$ for all $n \in \N$, and the
similar relations for $(V_{n*},V_n)$ and $(V_{n*}^{\prime},V_n^{\prime})$. Let $N \in \N$ be such that
$2^{-N} \leq \epsilon/2$.
By the above, we may
find a partition $\{I,J\}$ of $[N,+\infty)$,  vectors $u \in \sum_{n \in I}U_n^{\prime}$,
$v \in \sum_{n \in J}V_n^{\prime}$, and functionals $u^* \in \sum_{n \in I}U_{n*}^{\prime}$,
$v^* \in \sum_{n \in J}V_{n*}^{\prime}$, such that  $\norm{u-v}\norm{u^*-v^*}< \epsilon
|(u^*-v^*)(u+v)|$.

Let for $n \geq N$, $(W_{n*},W_n)=(U_{n*},U_n)$ if $n \in I$, and
$(W_{n*},W_n)=(V_{n*},V_n)$ if $n \in J$, and let $W_*=\sum_{n \geq N}W_{n*}$; let
$(W_{n*}^{\prime},W_n^{\prime})$ and $W_*^{\prime}$ be defined in a similar way. Let
also $X_N=[e_i, i \in \cup_{n \geq N}E_n]$.
We have therefore
$$\norm{u-v} < \epsilon \norm{u+v}_{W_*^{\prime}}.$$

Since $\sum_{n \geq N}\delta_n \leq \epsilon$, we find by Lemma \ref{perturbations} a map $T$
from $\sum_{n \geq N}W_n^{\prime}$ onto $\sum_{n \geq N}W_n$ such that
$T(W')=W$, 
$T(X_N \cap {W_{*}^{\prime}}^{\perp})=X_N \cap W_{*}^{\perp}$, and with
$\norm{T}\norm{T^{-1}} \leq (1+\epsilon)(1-\epsilon)^{-1} \leq 2$. 
If let $x=Tu \in \sum_{n \in I}U_n$ and $y=Tv \in \sum_{n \in J}V_n$, we have 

$$\norm{x-y} < 2\epsilon \norm{x+y}_{W_*}.$$

This means that we may pick some normalized functional $w^* \in W_*$, therefore $w^*=x^*-y^*$ with
$x^* \in \sum_{n \in I}U_{n*}$, $y^* \in \sum_{n \in J}V_{n*}$, such that

$$\norm{x-y} < 2\epsilon|(x^*-y^*)(x+y)|=2\epsilon |x^*(x)-y^*(y)|,$$

and we deduce that $A({\cal U},{\cal V}) \leq 2\epsilon$.
 \pff

\section{Remarks and open questions}

\paragraph{Remark 21} Let $Y$ be an FDD-block quotient of subspace of $X$.
 To check whether  $Y$ has the restricted QHI property, we have checked the formally stronger result
that the angle is zero between any two QS-pairs associated to FDD-block 
quotient of subspaces of Y with sequences of blocks having the same partition. 
We note that by our dichotomy theorem, these two notions are equivalent up to passing to
a quotient of subspace. Indeed if $X$
 is QHI restricted to block-subspaces, then no quotient of $X$ by an FDD-block subspace can contain an
unconditional basic sequence, and therefore $X$  must contain a quotient of subspace with the stronger "angle zero"
property.

\paragraph{Question 22} Is it possible to improve Theorem \ref{tadam} to suppress the restriction
to FDD-block subspaces? The restricting condition is not only technical. 
By a result of S. Argyros, A. Arvanitakis and A. Tolias \cite{AAT}, the distinction between
general quotient spaces and quotient by FDD-block subspaces can be essential: 
there exists a separable dual space $X$ with a Schauder basis, such that quotients with $w^*$-closed kernels are HI, yet every quotient
has a further quotient
isomorphic to $l_2$. Since FDD-block subspaces of $X$ are $w^*$-closed, this space has the restricted
QHI property, but it is not QHI by the $l_2$-saturation property.

Contrary to the case of subspaces, it does not seem clear that, in a space with a Schauder basis,
 QS-spaces may be approximated by FDD-block
quotient of subspaces, that is, that for any QS-space, there is a further QS-space, which is an arbitrary
small perturbation of an FDD-block quotient of subspace.

%
%

\paragraph{Remark 23} As was noticed in the introduction, HI spaces can fall in either side of the dichotomy
in Theorem \ref{theo1}.
The example of $X_{GM}$ is QHI, while the examples of \cite{AF} have an unconditional quotient. 
The dual $X_{uh}^*$ of the reflexive space $X_{uh}$ of 
  Argyros and Tolias \cite{AT2} has the following quite interesting mixed property.
 Any of its quotients  has
a further quotient with an unconditional basis, \cite{AT2} Proposition 3.6. On the other hand
it is HI, \cite{AT2} Proposition 5.11, and  it is saturated with
 QHI subspaces. This last fact was indicated to us by S. Argyros and the proof is as follows.
Consider any block subspace
of $X_{uh}^*$. Keeping a half of
the
vectors of the block basis, and denoting the space generated by them $Y$, we get
that the annihilator of any subspace $Z$ of $Y$ must contain an
infinite subsequence of the basis. Therefore \cite{AT2} Proposition 6.3. applies
to obtain that $X_{uh}/Z^{\perp}$ is HI. This means that every infinite dimensional quotient
of $X_{uh}/Y^{\perp}$ is HI, therefore $X_{uh}/Y^{\perp}$ is QHI. By reflexivity it follows that
$Y \simeq (X_{uh}/Y^{\perp})^*$ is QHI.

 \

We say that a Banach space $X$ is {\em unconditionally QS-saturated (resp. QS-saturated with HI subspaces)}
  if any infinite dimensional QS-space of $X$ has a further QS-space
 with an unconditional basis (resp. which is HI).

By Odell's result \cite{O}, if a space $X$ has a shrinking unconditional FDD, then every quotient of $X$ must be
unconditionally saturated, and therefore $X$ must be unconditionally QS-saturated.
It remains unknown whether there exists a HI space which is unconditionally QS-saturated.  Therefore we ask: 

\

Does every HI space contain a quotient of subspace which
is QHI? or which has the restricted QHI property?

\paragraph{Remark 24} Our dichotomy theorem, the result of Odell \cite{O}, and the
 remark after Proposition \ref{QHIreflexif}
 imply the following: if 
 $X$ is reflexive and QS-saturated with HI spaces, then some quotient
of subspace of $X$ is saturated with subspaces with HI dual.

In this direction, we recall the question of S. Argyros:

\

Does there exist a reflexive HI space $X$, such that no subspace of $X$ has a HI dual? 

\

\

\

Valentin Ferenczi,

 Equipe d'Analyse Fonctionnelle,

Universit\'e Pierre et Marie Curie - Paris 6,

Bo\^\i te 186, 4, Place Jussieu,

 75252, Paris Cedex 05,

 France.

\

 E-mail: ferenczi@ccr.jussieu.fr.

\end{document}